\documentclass[notitlepage]{amsart}

\usepackage{amsfonts}
\usepackage[latin1]{inputenc}
\newtheorem{pro}{Proposition}[section]
\newtheorem{teo}[pro]{Theorem}
\newtheorem{defi}[pro]{Definition}
\newtheorem{lem}[pro]{Lemma}
\newtheorem{cor}[pro]{Corollary}
\newtheorem{remark}[pro]{Remark}

\newcommand{\pd}{{\mathrm{pd}}}

\newcommand{\modu}{{\mathrm{mod}}}

\newcommand{\Ker}{{\mathrm{Ker}}}

\newcommand{\proj}{{\mathrm{proj}}}
\newcommand{\Proj}{{\mathrm{Proj}}}

\newcommand{\findim}{{\mathrm{fin.dim}}}

\newcommand{\gldim}{{\mathrm{gl.dim}}}
\newcommand{\Hom}{{\mathrm{Hom}}}

\newcommand{\Ext}{{\mathrm{Ext}}}
\newcommand{\Phidim}{{\Phi\,\mathrm{dim}}}
\newcommand{\Psidim}{{\Psi\,\mathrm{dim}}}

\newcommand{\T}{{\mathcal{T}}}
\newcommand{\C}{{\mathcal{C}}}
\newcommand{\A}{{\mathcal{A}}}

\newcommand{\add}{{\mathrm{add}}}

\newcommand{\rad}{{\mathrm{rad}}}

\newenvironment{dem}{\noindent {\bf Proof.}}{\hfill $\Box$\\}

\usepackage{latexsym,amssymb,amscd}
\usepackage{amsmath}
\usepackage[all,cmtip]{xy}
\usepackage[usenames]{color}
\usepackage{tikz-cd}
\usepackage{xypic}

\begin{document}

\title[Pullback diagrams, syzygy finite classes and Igusa-Todorov algebras]{ Pullback diagrams, syzygy finite classes and Igusa-Todorov algebras}
\author{Diego Bravo$^{1}$, Marcelo Lanzilotta$^{1}$ and Octavio Mendoza$^{2}$}
\footnote{The third author thanks the project PAPIIT-Universidad Nacional Aut\'onoma de M\'exico IN103317.}
 
\maketitle

\begin{center}
$^{1}$Instituto de Matem\'atica y Estad\'istica ``Rafael Laguardia'', Universidad de la Rep\'{u}blica, Uruguay.\\
$^{2}$Instituto de Matem\'aticas, Universidad Nacional Aut\'onoma de M\'exico, M\'exico.\\

dbravo@fing.edu.uy\\
marclan@fing.edu.uy\\
omendoza@matem.unam.mx
\end{center}

%\centerline{\small{We dedicate this paper to Eduardo Marcos on
 %his sixtieth birthday}}
 
\begin{abstract}
\noindent For an abelian category $\mathcal{A}$, we define the category PEx($\mathcal{A}$) of pullback diagrams of short exact sequences in 
$\mathcal{A}$, as a subcategory of the functor category Fun($\Delta, \mathcal{A}$) for a fixed diagram category $\Delta$. For any object
$M$ in ${\rm PEx}(\mathcal{A}),$ we prove the existence of a short exact sequence $0 {\to} K {\to} P {\to} M {\to} 0$ of functors, where 
the objects are in PEx($\mathcal{A}$) and $P(i) \in {\rm Proj(\mathcal{A})}$ for any $i \in \Delta$.  As an application, we prove that if 
$(\mathcal{C}, \mathcal{D}, \mathcal{E})$ is a triple of syzygy finite classes of objects in $\modu\,\Lambda$ satisfying some special conditions, then 
$\Lambda$ is an Igusa-Todorov algebra. Finally, we study lower triangular matrix Artin algebras and determine in terms of their components, under reasonable hypothesis, when these algebras are syzygy finite or Igusa-Todorov.
\ 


\end{abstract}

\section{Introduction}

We work on the setting of  Artin algebras and all the modules to be considered are finitely generated left modules. Given an Artin algebra 
$\Lambda,$ we  denote by $\modu\,\Lambda$ the category of finitely generated left $\Lambda$-modules.  The finitistic dimension of  $\Lambda$ is 
$\findim\,\Lambda:=\sup\{\pd\,M \hspace{0.1cm} | \hspace{0.1cm} M\in \mbox{mod}\,\Lambda \quad \mbox{and pd}\,M<\infty \}$, where 
$\pd\,M$ stands for the projective dimension of $M$. 
\

We recall that the finitistic dimension conjecture states that  $\findim\,\Lambda$ is finite, for any Artin algebra $\Lambda.$
 It is worth mentioning that the finitistic dimension conjecture is still open and it is one of the main problems in 
the representation theory of algebras. For more information about the history of the finitistic dimension conjecture, we suggest the reader to see 
in \cite{Z}. Until now, it is known that this conjecture is true for several  classes of algebras, among others: algebras with radical cube zero, 
monomial algebras, left serial algebras, weakly stably hereditary algebras and special biserial algebras. A large class of algebras, containing the 
mentioned before, is the class of Igusa-Todorov algebras. 
\

Igusa-Todorov algebras were introduced by J. Wei in \cite{W}. In this paper, J. Wei provides some methods to construct new classes of Igusa-Todorov algebras and thus we obtain algebras such that the finitistic dimension conjecture holds true. In the attempt to construct more Igusa-Todorov algebras, 
J. Wei introduces in \cite[Theorem 1.1]{W2} a new method that uses syzygy finite classes of modules.
\

 A class $\mathcal{C}$ of $\Lambda$-modules is said to be of finite representation type, if there is some 
$N\in\modu\,\Lambda$ such that  $\add\,\mathcal{C}=\add\,N,$ where $\add\,\mathcal{C}$ is the class of all $\Lambda$-modules which are direct 
summands of finite direct sums of objects in $\C.$ If $\mathcal{C}=\mbox{mod}\,\Lambda$, it is said that $\Lambda$ is of finite representation type. We denote by $\Omega^{i}$ the ith-syzygy operator. A class $\mathcal{C}$ of $\Lambda$-modules is said to be \textit{n-syzygy finite}, for some $n\geq 0,$ if the class 
 $\Omega^n\mathcal{C}$ is of finite representation type. In this case, 
the $\Lambda$-module $M:=\displaystyle\oplus_{i=1}^{t}M_i$ is called the \textit{$n$-syzygy representative} $\Lambda$-module for $\C,$ where $M_1, \ldots, M_t$ are the pairwise non-isomorphic indecomposable direct summands of $\Lambda$-modules in $\Omega^n\mathcal{C}$. If $\modu\,\Lambda$ is $n$-syzygy finite, we say that $\Lambda$ is $n$-syzygy finite. We recall that,  for the given classes $\C$ and 
$\mathcal{D}$ of $\Lambda$-modules,  $\C \ast \mathcal{D}$ is the class of all objects $M\in\modu\,\Lambda$ for which there is an exact sequence $0\to C\to M\to D\to 0$ with $C\in\C$ and $D\in\mathcal{D}.$
\

In this paper, we generalize the ideas proposed by J. Wei in the proof of \cite[Theorem 1.1]{W2}, where the syzygy finite classes play an important role.  The main result, going in this direction, is given in section 3 and it is the following one.
\

{\bf Theorem \ref{IT}} Let $\Lambda$ be an Artin algebra and let $n,p$ be non-negative integers. Let $\mathcal{C},$ 
$\mathcal{D},$ $\mathcal{E}$  and $\mathcal{K}$ be classes of objects in $\modu\,\Lambda$ such that
 $\C$ and $\mathcal{E}$ are $n$-syzygy finite, and 
$\mathcal{K}$ is $p$-syzygy finite. If the following  conditions hold
\begin{itemize}
\item[(TS1)]  $\;\Omega(\modu\,\Lambda)\subseteq\C \ast \mathcal{D};$ 
\item[(TS2)]  for any $D \in \mathcal{D}$ there is an exact sequence $0 \to K \to E \to D \to 0$, where $E \in \mathcal{E}$ and $K\in\mathcal{K};$ and
\item[(TS3)] there is some $j>0$ such that $\Ext^j_{\Lambda}(\mathcal{E},\Omega^{j-1}\C)=0;$
\end{itemize}
then $\Lambda$ is a $m$-Igusa-Todorov algebra, for $m:=\max\{p,n,j\}+1.$

The main ingredient in the proof of Theorem \ref{IT}, relies on the construction of an exact sequence of functors, for any pullback diagram of 
$\Lambda$-modules $A,$  of the form $0\to K\to P\to A\to 0.$ In this exact sequence, $K$ and $P$ are also pullback diagrams and $P$ has in any entry a projective $\Lambda$-module. The existence of such exact sequence, and all the details concerned with that, is discussed in section 2.
\

The study of lower triangular matrix Artin algebras is also included in the development of this paper. This class of algebras is the concern of section 4.
A lower triangular matrix Artin $R$-algebra $\Lambda$ is  of the form $\Lambda={\begin{pmatrix} T & 0 \\ {}_{U}M_T & U\end{pmatrix}},$ where $T$ and $U$ are Artin $R$-algebras,  
 $M$ is a $U$-$T$-bimodule, where $R$ acts centrally on $M,$ and furthermore $M$ is 
finitely generated over $R.$ The description of the $\Lambda$-modules can be given, see the details in \cite[III.2]{ARS}, in terms of 
triples $(A,B,f),$ where $A$ is a $T$-module, $B$ is a $U$-module and $f:M\otimes A\to B$ is a morphism of $U$-modules. For this triangular 
matrix algebra $\Lambda,$ it is very useful to consider the following classes of $\Lambda$-modules 
$$\mathcal{T}:=\{\overline{A}:=(A,0,0)\in 
 \modu\,\Lambda\;:\; A\in \modu\,T\},$$  $$\mathcal{U}:=\{\underline{B}:=(0,B,0)\in \modu\,\,\Lambda\;:\; B\in \modu\,U\}.$$
Note that, for any   $(A,B,f)\in\modu\,\Lambda,$ we have the exact sequence $0\to \underline{B}\to (A,B,f)\to\overline{A}\to 0,$ and thus 
$\modu\,\Lambda=\mathcal{U} *\T.$ Then, as an 
application of Theorem \ref{IT}, we obtain the following result.
\
 
{\bf Corollary \ref{LT1}}   Let $\Lambda={\begin{pmatrix} T & 0 \\ _{U}M_T & U\end{pmatrix}}$ be a triangular matrix Artin $R$-algebra such that $\T$ is $n$-syzygy finite and $\mathcal{U}$ is $m$-syzygy finite. Then, $\Lambda$ is a $k$-Igusa-Todorov algebra for $k:=\max\{m,n\}+1.$
\

In order to check the $n$-syzygy finiteness condition, for the classes $\mathcal{U}$ and $\mathcal{T},$ 
we need to have a description of $\Omega_\Lambda$ in terms of $\Omega_T$ and $\Omega_U.$ In general, it is a hard problem to attemp. However, 
under some natural hypothesis on the bimodule $M,$ it is possible to handle this situation and thus we get a very nice description of $\Omega_\Lambda.$ 
To state this description, we consider the following notation. For any Artin algebra $\Gamma$ and any $A\in\modu\,\Gamma,$ the exact sequence 
$$\cdots\to P^A_n\to P^A_{n-1}\to\cdots \to P^A_1\to P^A_0\to A\to 0$$ 
denotes the minimal projective resolution of $A.$ Thus, for any 
 $k\geq 1,$ we have the exact sequence $0\to\Omega^kA\xrightarrow{i^A_k} P^A_{k-1}\xrightarrow{\pi^A_{k-1}} \Omega^{k-1}A\to 0$ of $\Gamma$-modules.
 
{\bf Lemma \ref{triang}}
Let $\Lambda={\begin{pmatrix} T & 0 \\ _{U}M_T & U\end{pmatrix}}$ be a triangular matrix Artin $R$-algebra such that ${}_U M\in\proj\,(U)$
 and $M_T \in\proj\,(T).$ Then, the $n$-th syzygy of 
$(A, B, f)\in\modu\,\Lambda$ is
$$\Omega^n(A, B, f)=(\Omega^nA, M\otimes P^A_{n-1}, M\otimes i^A_{n}) \oplus (0, \Omega^nB,0).$$
The above lemma  plays a  fundamental role in the proof of the following pair of results obtained in this paper.

{\bf Theorem \ref{teotriang}}   Let $\Lambda={\begin{pmatrix} T & 0 \\ {}_{U}M_T & U\end{pmatrix}}$ be a triangular matrix Artin $R$-algebra such that ${}_U M\in\proj\,(U)$ and $M_T \in\proj\,(T).$ Then, 
for a positive integer $n,$ the following statements are 
equivalent.
\begin{itemize}
\item[(a)] The algebras $T$ and $U$ are $n$-syzygy finite.
\item[(b)] The classes $\mathcal{T}$ and $\mathcal{U}$ are $n$-syzygy finite.
\item[(c)] The algebra $\Lambda$ is $n$-syzygy finite.
\end{itemize}
Moreover, if one of the above equivalent conditions hold, then $\Lambda$ is a $(n-1)$-Igusa-Todorov algebra.
\

{\bf Thoerem \ref{nITtr}} Let $\Lambda={\begin{pmatrix} T & 0 \\ _{U}M_T & U\end{pmatrix}}$ be a triangular 
matrix Artin $R$-algebra such that ${}_U M\in\proj\,(U)$ and $M_T \in\proj\,(T).$ Then
$\Lambda$ is $n$-Igusa-Todorov if, and only, if $T$ and $U$ are $n$-Igusa-Todorov. 
\

One situation where the needed hypothesis, in the last two theorems, hold true is in the particular case that $\Gamma:=T=U=M.$ Thus, $\Lambda$ 
becomes into a $2$-triangular matrix algebra $T_2(\Gamma).$ We can generalize this notion as follows.
\

For an Artin $R$-algebra $\Gamma,$   the $k$-triangular matrix $R$-algebra 
$T_k(\Gamma)$ is the matrix of 
size $k\times k$ with entries $[T_k(\Gamma)]_{i,j}:=\Gamma$ for $i\geq j,$ and 
$[T_k(\Gamma)]_{i,j}:=0$ otherwise. By using induction, the two theorems above allow us to prove the following result, and the last one, of this paper.

{\bf Corollary} For an Artin $R$-algebra $\Gamma,$ the following statements are equivalent.
\begin{itemize}
\item[(a)] $\Gamma$ is $n$-syzygy finite (resp. $n$-Igusa-Todorov).
\item[(b)] $T_k(\Gamma)$ is $n$-syzygy finite (resp. $n$-Igusa-Todorov), for any $k\geq 2.$
\item[(c)] $T_k(\Gamma)$ is $n$-syzygy finite (resp. $n$-Igusa-Todorov), for some $k\geq 2.$
\end{itemize}

\section{Pullback diagrams and functor categories}

In this section we introduce the category of pullback diagrams of exact sequences in an abelian category $\A.$  For any pullback diagram $A$ of 
objects in $\A,$  we prove the existence of an 
exact sequence of functors $0\to K\to P\to A\to 0,$ where $K$ and $P$ are also pullback diagrams and $P$ has in any entry a projective object in 
$\A.$ We denote by $\Proj\,(\A)$ the class of all projective objects in the abelian category $\A.$
\
 
In order to introduce the category of pullback diagrams in $\A,$ we start  by fixing the following very 
 special  quiver $Q$
$$Q:\qquad\qquad\begin{CD}
     @. e_7 @>\alpha_{78}>> e_8 \\
   @.    @V\alpha_{75}VV  @VV\alpha_{86}V     \\
 e_4 @>\alpha_{45}>> e_5 @>\alpha_{56}>> e_6 \\
   @V\alpha_{41}VV    @V\alpha_{52}VV  @VV\alpha_{63}V \\
e_1 @>>\alpha_{12}> e_2 @>>\alpha_{23}> e_3 \\
 \end{CD}$$
Consider the $\mathbb{Z}$-category $\Delta:=\mathbb{Z}Q/\sim$, where $\sim$ is given by the commutativity relations of the squares in the diagram $Q$. Recall that Ob($\Delta)=\{e_1, \ldots, e_8 \}$ and $\Delta(e_i, e_j)$ is the free abelian group generated by the paths in $Q$ from $e_i$ to $e_j$, modulo the commutativity relations in the diagram above.
\

Let $\mathcal{A}$ be an abelian category, and let $\mathrm{Fun}(\Delta, \mathcal{A})$ be the category of all additive covariant functors from 
$\Delta$ to $\mathcal{A}.$ Since $\mathcal{A}$ is abelian, it follows that $\mathrm{Fun}(\Delta, \mathcal{A})$ is also an abelian category. 

\begin{defi}
Let $\mathcal{A}$ be an abelian category. The  category of {\bf pullback diagrams} of exact sequences in $\A,$ denoted by $\mathrm{PEx}(\mathcal{A}),$ is the full subcategory of $\mathrm{Fun}(\Delta, \mathcal{A})$ whose objects are the functors $M:\Delta\to\A$ satisfying the following conditions

\begin{itemize}
	\item $M(\alpha_{41})$ and  $M(\alpha_{78})$ are isomorphisms;
	\item $M(\alpha_{12}), M(\alpha_{45}), M(\alpha_{75}),$ and  $M(\alpha_{86})$ are monomorphisms in 
	$\mathcal{A};$
	\item $M(\alpha_{52}), M(\alpha_{63}), M(\alpha_{56}),$ and $M(\alpha_{23})$ are epimorphisms in 
	$\mathcal{A};$ 
	\item $M$ sends the diagram $Q$ to an exact diagram in $\mathcal{A}$ (i.e with exact rows and columns).
	\end{itemize} 
	
\end{defi}

\begin{remark} Note that $\mathrm{PEx}(\mathcal{A})$ is closed under isomorphisms in $\mathrm{Fun}(\Delta, \mathcal{A}).$ However, in general, it is not closed under extensions.
\end{remark}

Indeed, let $k$ be a field. Consider the short exact sequence $0 {\to} K {\to} M {\to} A {\to} 0$ in ${\rm Fun}(\Delta, {\rm mod}\,k)$ with $K, A$ in ${\rm PEx}({\rm mod}\,k)$ described by the following commutative diagram in $\modu\,k$ 

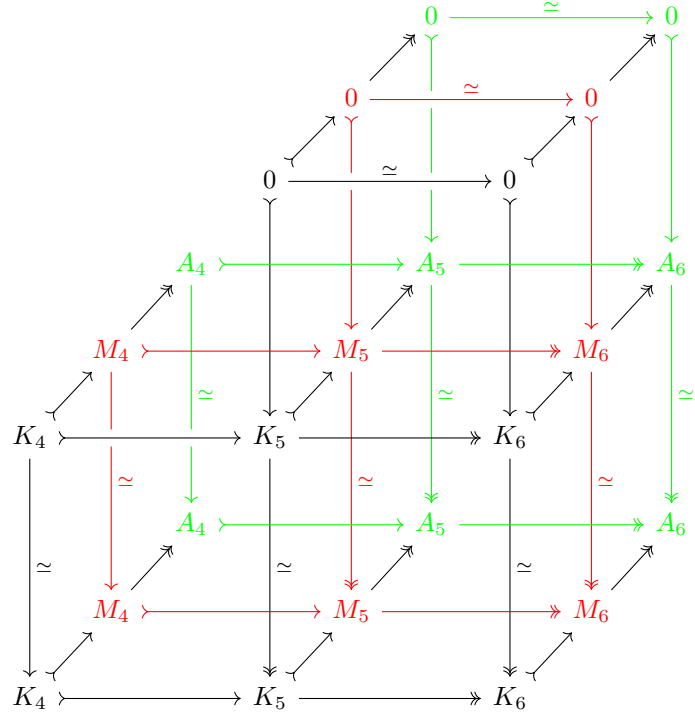
\[
\begin{tikzcd}[row sep=4ex, column sep=2ex]
& & &  & & & \textcolor{green}{0} \arrow[rrr, "\simeq", green] \arrow[ddd, rightarrowtail, green] &&& \textcolor{green}{0} \arrow[ddd, rightarrowtail, green]\\
&  & & & & \textcolor{red}{0} \arrow[rrr,red, "\simeq", crossing over] \arrow[ur,twoheadrightarrow] \arrow[ddd, rightarrowtail, red] &&& \textcolor{red}{0} \arrow[ur, twoheadrightarrow] \arrow[ddd, rightarrowtail, red] &\\
& & && 0 \arrow[rrr, "\simeq", crossing over] \arrow[ur, rightarrowtail] \arrow[ddd, crossing over, rightarrowtail] &&& 0 \arrow[ur, rightarrowtail] \arrow[ddd, rightarrowtail] &&\\
&&  & \textcolor{green}{A_4} \arrow[rrr, rightarrowtail, green] \arrow[ddd, "\simeq", green] & && \textcolor{green}{A_5} \arrow[rrr, twoheadrightarrow, green] \arrow[ddd, twoheadrightarrow, green, "\simeq"] &&& \textcolor{green}{A_6} \arrow[ddd, twoheadrightarrow, green, "\simeq"]\\
& &  \textcolor{red}{M_4} \arrow[ur, twoheadrightarrow] \arrow[ddd, "\simeq", red] \arrow[rrr, rightarrowtail, red] & && \textcolor{red}{M_5} \arrow[rrr, twoheadrightarrow, red] \arrow[ur, twoheadrightarrow] \arrow[ddd, twoheadrightarrow, red, "\simeq"] &&& \textcolor{red}{M_6} \arrow[ur, twoheadrightarrow] \arrow[ddd, twoheadrightarrow, red, "\simeq"] & &\\
& K_4 \arrow[rrr, rightarrowtail, crossing over] \arrow[ur, rightarrowtail] \arrow[ddd, "\simeq"] & && K_5 \arrow[rrr, twoheadrightarrow] \arrow[ur, rightarrowtail] \arrow[ddd, crossing over, twoheadrightarrow, "\simeq"] &&& K_6 \arrow[ur, rightarrowtail] \arrow[ddd, twoheadrightarrow, crossing over, "\simeq"]& &\\
&   &  & \textcolor{green}{A_4} \arrow[rrr, rightarrowtail, green] & && \textcolor{green}{A_5} \arrow[rrr, twoheadrightarrow, green] &&& \textcolor{green}{A_6}\\
&  & \textcolor{red}{M_4} \arrow[ur, twoheadrightarrow] \arrow[rrr, rightarrowtail, red] & && \textcolor{red}{M_5} \arrow[ur, twoheadrightarrow] \arrow[rrr, twoheadrightarrow, red] &&& \textcolor{red}{M_6} \arrow[ur, twoheadrightarrow] & & \\
& K_4 \arrow[rrr, rightarrowtail] \arrow[ur, rightarrowtail]   & && K_5 \arrow[rrr, twoheadrightarrow] \arrow[ur, rightarrowtail] &&& K_6 \arrow[ur, rightarrowtail] &  & \\
\end{tikzcd}
\]
Since the objects are finite dimensional $k$-vector spaces, the morphisms of the mid level diagram can 
be written as follows
\[
\begin{tikzcd}[ampersand replacement=\&]
K_4 \arrow[rr, rightarrowtail, "{\begin{pmatrix} 1\\ 0\end{pmatrix}}"] \arrow[ddd, rightarrowtail, "{\begin{pmatrix} 1\\ 0\end{pmatrix}}"] \&\& K_4 \oplus K_6 \arrow[rr, twoheadrightarrow, "{\begin{pmatrix} 0 & 1\end{pmatrix}}"] \arrow[ddd, rightarrowtail, "{\begin{pmatrix} 1 & 0 \\ 0 & 0\\ 0 & 1\\ 0 & 0\end{pmatrix}}"] \&\& K_6 \arrow[ddd, rightarrowtail,"{\begin{pmatrix} 1\\ 0\end{pmatrix}}"]\\
\&\& \&\&\\
\&\& \&\&\\
K_4 \oplus A_4 \ar[rr, rightarrowtail, "{\begin{pmatrix} 1 & x_1 \\ 0 & 1\\ 0 & x_2\\ 0 & 0\end{pmatrix}}"] \arrow[ddd, twoheadrightarrow, "{\begin{pmatrix} 0 & 1\end{pmatrix}}"]\&\& K_4 \oplus A_4 \oplus K_6 \oplus A_6 \arrow[rr, twoheadrightarrow, "{\begin{pmatrix} 0 & y_1 & 1 & y_2 \\ 0 & 0 & 0 & 1\end{pmatrix}}"] \arrow[ddd, twoheadrightarrow, "{\begin{pmatrix} 0 & 1 & 0 & 0 \\ 0 & 0 & 0 & 1\end{pmatrix}}"] \&\& K_6 \oplus A_6 \arrow[ddd, twoheadrightarrow, "{\begin{pmatrix} 0 & 1\end{pmatrix}}"] \\
\&\& \&\&\\
\&\& \&\&\\
A_4 \arrow[rr, rightarrowtail, "{\begin{pmatrix} 1\\ 0\end{pmatrix}}"] \&\& A_4 \oplus A_6 \arrow[rr, twoheadrightarrow, "{\begin{pmatrix} 0 & 1\end{pmatrix}}"] \&\& A_6\\
\end{tikzcd}
\]
This diagram commutes and thus the big diagram commutes, but in general, the middle row is not a short exact sequence. For instance, take $A_4\neq 0,$ $K_6\neq 0,$ $x_2=0$ and $y_1\neq 0.$ Thus, we have shown the existence of a short exact sequence $0 {\to} K {\to} M {\to} A {\to} 0$ in ${\rm Fun}(\Delta, {\rm mod}\,k)$ with $K, A$ in ${\rm PEx}({\rm mod}\,k)$, but $M \notin {\rm PEx}({\rm mod}\,k).$

\begin{lem}\label{sum} The class
	$\mathrm{PEx}(\mathcal{A})$ is closed under finite direct sums and direct summands in  
	$\mathrm{Fun}(\Delta,\mathcal{A}).$
\end{lem}
\begin{dem} It follows from the following fact. A sequence 
\[ A_1 \bigoplus A_2\xrightarrow{\begin{pmatrix} f_1 & 0 \\ 0 & f_2 \end{pmatrix}} 
B_1 \bigoplus B_2\xrightarrow{ \begin{pmatrix} g_1 & 0 \\ 0 & g_2 \end{pmatrix} }C_1 \bigoplus C_2   \]
in the abelian category $\A$ is exact, if and only if, the sequences
$A_1\xrightarrow{f_1}B_1\xrightarrow{g_1}C_1$ and $A_2\xrightarrow{f_2}B_2\xrightarrow{g_2}C_2$ are exact.
\end{dem}

\begin{defi} Let $\A$ be an abelian category. We denote by ${\rm PEx}({\rm Proj}\,(\mathcal{A}))$ the full subcategory of ${\rm PEx}(\mathcal{A})$ whose objects are the functors $P:\Delta\to\A$ such that $P(e)\in\Proj\,(\A),$ for any object $e \in \Delta$. We also consider $\Proj\,(\mathrm{PEx}(\A)),$ the class of all projective objects in 
the category $\mathrm{PEx}(\A).$
\end{defi}

\begin{cor}\label{sumproj} The class ${\rm PEx}({\rm Proj}\,(\mathcal{A}))$ is closed under finite direct sums 
and direct summands in  
	$\mathrm{Fun}(\Delta,\mathcal{A}).$
\end{cor}
\begin{dem} It follows from Lemma \ref{sum} and the fact that ${\rm Proj}\,(\mathcal{A})$ is closed under finite direct sums 
and direct summands in  $\A.$
\end{dem}

The following result is the main ingredient in the proof of the Theorem \ref{IT}.

\begin{pro}\label{ses}
Let $\mathcal{A}$ be an abelian category with enough projectives. Then, for any $A\in {\rm PEx}(\mathcal{A})$ there is a short exact sequence $0 \to K \to P \to A \to 0$ in 
$\mathrm{Fun}(\Delta, \A)$ with objects in $\mathrm{PEx}(\A)$ and 
$P \in \mathrm{PEx}(\Proj\,(\A)).$	
\end{pro}
\begin{dem} Let $A\in {\rm PEx}(\mathcal{A}).$ Thus $A(Q)$ is the following exact and commutative diagram in $\A$
	$$\begin{CD}
	@. @. 0 @. 0 @. @.\\
	@.   @.    @VVV  @VVV     @.\\
	@.   @. A_7 @>{\simeq}>> A_8 @.\\
	@.   @.    @VVV  @VVV     @.\\
	0 @>>> A_4 @>>> A_5 @>>> A_6 @>>> 0\\
	@.   @VV{\simeq}V    @VVV  @VVgV     @.\\
	0 @>>> A_1 @>>> A_2 @>>> A_3 @>>> 0\\
	@.   @.    @VVV  @VVV     @.\\
	@. @. 0 @. 0. @. @.
	\end{CD}$$
Therefore, $A(Q)$ can be seen as the pullback diagram for the short exact sequence $$0 \longrightarrow A_1 \longrightarrow A_2 \longrightarrow A_3 \longrightarrow 0$$	and the epimorphism $g: A_6 \to A_3$ in $\mathcal{A}$. Since $\mathcal{A}$ has enough projectives, there are epimorphisms $p_i: P_i \to A_i$ for $i=3,4,7.$ Let  $P_8:=P_7,$ $P_1:=P_4,$ $P_{ij}:=P_i \oplus P_j$ and $P_{ijk}:=P_i \oplus P_j \oplus P_k$. Thus, we can construct the following commutative diagram 

\[
\begin{tikzcd}[row sep=4ex, column sep=4ex]
  & & & A_7 \arrow[rr, "\simeq"] \arrow[dd, rightarrowtail] && A_8 \arrow[dd, rightarrowtail] & \\
   & & P_7 \arrow[rr, equal, near end,crossing over] \arrow[ur, twoheadrightarrow] \arrow[dd, crossing over, rightarrowtail] && P_8 \arrow[ur, twoheadrightarrow] \arrow[dd, rightarrowtail, crossing over] && \\
 & A_4 \arrow[rr, near start, "w", rightarrowtail] \arrow[dd,"\simeq", near end]  && A_5 \arrow[rr, near start, "\alpha", twoheadrightarrow] \arrow[dd, near start, "\beta", twoheadrightarrow] && A_6 \arrow[dd, near end, "g", twoheadrightarrow]\\
 P_4 \arrow[ur, "p_4", twoheadrightarrow] \arrow[dd, near end, crossing over, equal] \arrow[rr, near end, "r", rightarrowtail] & & P_{438} \arrow[rr, near end, "t", , twoheadrightarrow] \arrow[ur, dashrightarrow, "f", twoheadrightarrow] \arrow[dd, near start, "v", , twoheadrightarrow, crossing over] && P_{38} \arrow[ur, "s", , twoheadrightarrow] \arrow[dd,crossing over, twoheadrightarrow] & & \\
    & A_1 \arrow[rr, rightarrowtail] & & A_2 \arrow[rr, near start, "h", , twoheadrightarrow] && A_3\\
  P_1 \arrow[ur, , twoheadrightarrow] \arrow[rr, rightarrowtail] & & P_{13} \arrow[ur, "u", twoheadrightarrow] \arrow[rr, , twoheadrightarrow] && P_3 \arrow[ur, twoheadrightarrow] & & \\
\end{tikzcd}
\]
where all rows and columns are short exact sequences. We show the existence of an epimorphism 
$f: P_{438} \to A_5$ which makes the entire diagram commute. Indeed, since $huv = gst,$ $A_5$ is a pullback, and $p_4$ and $s$ are epimorphisms, there exists an epimorphism $f: P_{438} \to A_5$ such that $st=\alpha f$ and $uv=\beta f$. It remains to show the commutativity of two squares with the map $f$ as a side.
\

 Using that $A_5$ is a pullback and $gstr=huvr$, it follows the existence of a unique map $f': P_4 \longrightarrow A_5$ such that $\alpha f'=str$ and $\beta f'=uvr$. Since $\alpha fr=str=\alpha wp_4$ and $\beta fr=uvr=\beta wp_4$, we get $fr=wp_4$. The commutativity of the last square can be shown in a similar way. We have proven the existence of an epimorphism $P \twoheadrightarrow A$, where $P\in {\rm PEx}({\rm Proj}(\mathcal{A}))$. Note that $P(e_i)=P_i,$ for $i=1,3,4,7,$ $P(e_2)=P_{13},$ $P(e_5)=P_{438}$ and $P(e_6)=P_{38}.$

\noindent Consider now the syzygy $\Omega{A_i}:=\Ker\,(P(e_i)\twoheadrightarrow A_i).$  
Then, we can construct the following exact and commutative diagram  
\[
\begin{tikzcd}[row sep=4ex, column sep=2ex]
 & & &  & & & \textcolor{green}{A_7} \arrow[rrr, "\simeq", green] \arrow[ddd, rightarrowtail, green] &&& \textcolor{green}{A_8} \arrow[ddd, rightarrowtail, green]\\
 &  & & & & \textcolor{red}{P_7} \arrow[rrr,red, equal, crossing over] \arrow[ur,twoheadrightarrow] \arrow[ddd, rightarrowtail, red] &&& \textcolor{red}{P_8} \arrow[ur, twoheadrightarrow] \arrow[ddd, rightarrowtail, red] &\\
 & & && {\Omega A_7} \arrow[rrr, "\simeq", crossing over] \arrow[ur, rightarrowtail] \arrow[ddd, crossing over, rightarrowtail] &&& {\Omega A_8} \arrow[ur, rightarrowtail] \arrow[ddd, rightarrowtail] &&\\
 &&  & \textcolor{green}{A_4} \arrow[rrr, rightarrowtail, green] \arrow[ddd, "\simeq", green] & && \textcolor{green}{A_5} \arrow[rrr, twoheadrightarrow, green] \arrow[ddd, twoheadrightarrow, green] &&& \textcolor{green}{A_6} \arrow[ddd, twoheadrightarrow, green]\\
 & &  \textcolor{red}{P_4} \arrow[ur, twoheadrightarrow] \arrow[ddd, equal, red] \arrow[rrr, rightarrowtail, red] & && \textcolor{red}{P_{438}} \arrow[rrr, twoheadrightarrow, red] \arrow[ur, twoheadrightarrow] \arrow[ddd, twoheadrightarrow, red] &&& \textcolor{red}{P_{38}} \arrow[ur, twoheadrightarrow] \arrow[ddd, twoheadrightarrow, red] & &\\
 & {\Omega A_4} \arrow[rrr, rightarrowtail, crossing over] \arrow[ur, rightarrowtail] \arrow[ddd, "\simeq"] & && {\Omega A_5} \arrow[rrr, twoheadrightarrow] \arrow[ur, rightarrowtail] \arrow[ddd, crossing over, twoheadrightarrow] &&& {\Omega A_6} \arrow[ur, rightarrowtail] \arrow[ddd, twoheadrightarrow, crossing over]& &\\
 &   &  & \textcolor{green}{A_1} \arrow[rrr, rightarrowtail, green] & && \textcolor{green}{A_2} \arrow[rrr, twoheadrightarrow, green] &&& \textcolor{green}{A_3}\\
 &  & \textcolor{red}{P_1} \arrow[ur, twoheadrightarrow] \arrow[rrr, rightarrowtail, red] & && \textcolor{red}{P_{13}} \arrow[ur, twoheadrightarrow] \arrow[rrr, twoheadrightarrow, red] &&& \textcolor{red}{P_3} \arrow[ur, twoheadrightarrow] & & \\
 & {\Omega A_1} \arrow[rrr, rightarrowtail] \arrow[ur, rightarrowtail]   & && {\Omega A_2} \arrow[rrr, twoheadrightarrow] \arrow[ur, rightarrowtail] &&& {\Omega A_3} \arrow[ur, rightarrowtail] &  & \\
\end{tikzcd}
\]

 It is clear that the syzygy diagram commutes and thus the entire diagram commutes. We have shown the existence of $K\in {\rm PEx}(\mathcal{A})$, where $K(Q)$ is the above syzygy diagram, such that 
$0 \to K \to P \to A \to 0$ is a short exact sequence  in ${\rm Fun}(\Delta, \mathcal{A})$ with objects in ${\rm PEx}(\mathcal{A})$ and $P\in \mathrm{PEx}(\Proj\,(\A)).$ 
\end{dem}

\begin{cor}
Let $\A$ be an abelian category with enough projectives. Then 
$$\Proj\,({\rm PEx}(\A)) \subseteq {\rm PEx}(\Proj\,(\A)).$$	
\end{cor}
\begin{dem} It follows from Corollary \ref{sumproj} and Proposition \ref{ses}.
\end{dem}

\section{Igusa-Todorov algebras}

By studying the finitistic dimension conjecture, K. Igusa and G. Todorov  introduced two functions $\Phi$ and $\Psi$ from $\modu\,\Lambda$ to the 
non-negative integers \cite{IT}. These functions are nowadays known as the Igusa-Todorov functions, or IT-functions for short. For a further development 
of IT-functions, we recommend the reader to see in \cite{FLM, LM}.
 
 In the following lemma, we collect some basic well known properties of these functions.

\begin{lem}\label{desigualdades}
Let  $\Lambda$ be an Artin algebra. Then, the IT-functions $\Phi$ and $\Psi$ satisfy the 
following properties.
	
	\begin{itemize}
		\item[(a)] $\Psi(M)=\Phi(M)=\mbox{\emph{pd}}\,M,$ if $\pd\,M$ is finite.
		\item[(b)] $\Psi(M)=\Psi(N)$ and $\Phi(M)=\Phi(N)$ if $\add\,M=\add\,N.$
		\item[(c)] $\Psi(X) \leq \Psi(X\oplus Y)$  and $\Phi(X) \leq \Phi(X\oplus Y).$
		\item[(d)] $\Psi(X) \leq \Psi(\Omega^{n}X) + n$ and $\Phi(X) \leq \Phi(\Omega^{n}X) + n,$ for any $n\geq 0.$
		
		\item[(e)] For an exact sequence $0 \to X \to Y \to Z \to 0$  in $\modu\,\Lambda,$ we have that
		\begin{itemize}
		\item[($\mathrm{e_1}$)] if $\pd\,Z < \infty$ then $\pd\,Z \leq \Psi(X \oplus Y) + 1;$
		\item[($\mathrm{e_2}$)] if $\pd\,X < \infty$ then $\pd\,X \leq \Psi(\Omega Y \oplus \Omega Z) + 1;$
		\item[($\mathrm{e_3}$)] if $\pd\,Y < \infty$ then $\pd\,Y \leq \Psi(X \oplus \Omega\,Z) + 1.$
		\end{itemize}	
	\end{itemize}
\end{lem}
\begin{dem} See the proofs given in \cite{HLM1,IT,LM, Wa}.
\end{dem}

Related with the IT-functions, there are the $\Phi$-dimension and the $\Psi$-dimension which were introduced in \cite{HLM2} and  defined as follows
$$\Phidim\,\Lambda:=\sup\{\Phi(M)\;|\;M \in \mbox{mod}\,\Lambda\}\quad\mathrm{and}\quad\Psidim\,\Lambda:=\sup\{\Psi(M)\;|\;M \in \mbox{mod}\,\Lambda\}.$$

\noindent For an Artin algebra $\Lambda$, we have $\findim\,\Lambda \leq \Phidim\,\Lambda \leq \Psidim\,\Lambda \leq \gldim\,\Lambda$ where $\gldim\,\Lambda$ is the global dimension of 
$\Lambda.$

The next lemma is somehow well known and states that, for $n$-syzygy finite algebras, the $\Phi$-dimension is finite. Hence, by the above inequalities, the finitistic dimension of any $n$-syzygy finite algebra is finite

\begin{lem}\label{phifinite}
Let $\Lambda$ be a $n$-syzygy finite Artin algebra, and let $M$ be its $n$-syzygy representative $\Lambda$-module. Then $\Phidim\,\Lambda\leq \Phi(M)+n<\infty.$ 
\end{lem}
\begin{dem} Let $X \in \modu\,\Lambda$ and $M=\bigoplus_{i=1}^{t}M_i$, where $M_1, \ldots, M_t$ are all the pairwise non-isomorphic indecomposable direct summands of objects in 
$\Omega^n(\mbox{mod}\,\Lambda).$ 
Since $\Omega^{n}X = \bigoplus_{i=1}^{t}M_{i}^{\alpha_i},$ by Lemma \ref{desigualdades} we get
$$\Phi(X) \leq \Phi(\Omega^{n}X) + n = \Phi(\bigoplus_{i=1}^{t}M_{i}^{\alpha_i}) + n \leq \Phi(\bigoplus_{i=1}^{t}M_{i}) + n$$
and thus the lemma follows.
\end{dem}

Let $\Lambda$ be an Artin algebra and $n$ be a non-negative integer. Following \cite{W}, we say that 
$\Lambda$ is a \textbf{$n$-Igusa-Todorov algebra} if there exists $V\in\modu\,\Lambda$ such that for every $X \in \modu\,\Lambda$  there is a short exact sequence $$0 \longrightarrow V_{1} \longrightarrow V_{0} \longrightarrow \Omega^{n}X \oplus P \longrightarrow 0, $$
	
\noindent with $V_{1}, V_{0} \in \add\,V$ and $P$ is a projective $\Lambda$-module. In this case we say 
that $V$ is a $n$-Igusa-Todorov $\Lambda$-module.

\begin{lem}\label{ComIT} Let $\Lambda$ be a $n$-Igusa-Todorov algebra, and let $V$ be a $n$-Igusa-Todorov $\Lambda$-module. Then, the following conditions hold true.
\begin{itemize}
\item[(a)] \cite[Theorem 2.3]{W} $\findim\,\Lambda\leq \Psi(V)+n+1<\infty.$
\item[(b)] For all $X \in \modu\,\Lambda$  there is an exact sequence $0\to V_1\to V_0\to \Omega^n X\to 0,$ where $V_{1}, V_{0} \in \add\,V.$
\end{itemize}
\end{lem}
\begin{dem} (b) Let $0\to V_1\to V'_0\to \Omega^n X\oplus P\to 0$ be an exact sequence in $\modu\,\Lambda,$ where $V_{1}, V'_{0} \in \add\,V$ and $P\in\proj\,(\Lambda).$ Then, we have the following pullback diagram
\[\xymatrix{&& 0\ar[d]&0\ar[d]&\\
0\ar[r] & V_1 \ar@{=}[d]\ar[r] & V_0\ar[r]\ar[d]&\Omega^n X\ar[r]\ar[d]& 0\\
0\ar[r] & V_1\ar[r] & V'_0\ar[d] \ar[r]&\Omega^n X\oplus P\ar[d]\ar[r]& 0\\
&& P\ar@{=}[r]\ar[d]&P\ar[d]&\\
&& 0&0.&
}\]
Since $P\in\proj\,(\Lambda),$ the first column splits and thus $V_0\in\add\,V.$ Then the first row is the desired exact sequence.
\end{dem}

 In the next theorem, we  provide sufficient conditions, in terms of special subcategories 
of $\modu\,\Lambda,$ for an algebra $\Lambda$ to be $n$-Igusa-Todorov. For the given classes $\C$ and 
$\mathcal{D}$ of $\Lambda$-modules, we denote by $\C \ast \mathcal{D}$ the class of objects $M\in\modu\,\Lambda$ for which there is an exact sequence $0\to C\to M\to D\to 0$ with $C\in\C$ and $D\in\mathcal{D}.$

\begin{teo}\label{IT}
Let $\Lambda$ be an Artin algebra and let $n,p$ be non-negative integers. Let $\mathcal{C},$ 
$\mathcal{D},$ $\mathcal{E}$  and $\mathcal{K}$ be classes in $\modu\,\Lambda$ such that
 $\C$ and $\mathcal{E}$ are $n$-syzygy finite, and let
$\mathcal{K}$ be $p$-syzygy finite. If the following  conditions hold 
\begin{itemize}
\item[(TS1)]  $\;\Omega(\modu\,\Lambda)\subseteq\C \ast \mathcal{D};$ 
\item[(TS2)]  for any $D \in \mathcal{D}$ there is an exact sequence $0 \to K \to E \to D \to 0$, where $E \in \mathcal{E}$ and $K\in\mathcal{K};$ and
\item[(TS3)] there is some $j>0$ such that $\Ext^j_{\Lambda}(\mathcal{E},\Omega^{j-1}\C)=0;$
\end{itemize}
then $\Lambda$ is a $m$-Igusa-Todorov algebra, for $m:=\max\{p,n,j\}+1$.
\end{teo}
\begin{dem}
Let $X\in\modu\,\Lambda.$ By (TS1), there is an exact sequence $0\to C\to\Omega X\to D\to 0,$ with 
$C\in \mathcal{C}$ and $D\in \mathcal{D}.$ Moreover, from (TS2), we get an exact sequence 
$0\to K\to E\to D\to 0,$ where $E\in\mathcal{E}$ and $K\in\mathcal{K}.$
\

Consider the pull-back diagram 
$$\begin{CD}
@. @. 0 @. 0 @. @.\\
@.   @.    @VVV  @VVV     @.\\
@.   @. K @= K @.\\
@.   @.    @VVV  @VVV     @.\\
0 @>>> C @>>> G @>>> E @>>> 0\\
@.   @|    @VVV  @VVV     @.\\
0 @>>> C @>>> \Omega X @>>> D @>>> 0\\
@.   @.    @VVV  @VVV     @.\\
@. @. 0 @. 0. @. @.
\end{CD}$$ 
By Proposition \ref{ses}, we get the following exact and commutative diagram
$$\begin{CD}
@. @. 0 @. 0 @. @.\\
@.   @.    @VVV  @VVV     @.\\
@.   @. \Omega^{j-1}K @= \Omega^{j-1}K @.\\
@.   @.    @VVV  @VVV     @.\\
0 @>>> \Omega^{j-1}C @>>> \Omega^{j-1}G\oplus Q_2 @>>> \Omega^{j-1}E\oplus Q_3 @>>> 0\\
@.   @|    @VVV  @VVV     @.\\
0 @>>> \Omega^{j-1}C @>>> \Omega^{j}X\oplus Q_1 @>>> \Omega^{j-1}D @>>> 0\\
@.   @.    @VVV  @VVV     @.\\
@. @. 0 @. 0 @. @.
\end{CD}$$

\noindent where $Q_i$ is a projective $\Lambda$-module. By (TS3)  it follows that the middle row splits. Thus, the middle column give us the short exact sequence $$0 \longrightarrow \Omega^{j-1}K \longrightarrow \Omega^{j-1}C\oplus \Omega^{j-1}E\oplus Q_3 \longrightarrow \Omega^{j}X\oplus Q_1 \longrightarrow 0.$$
Since $m\geq j+1,$ from the preceding sequence, we get the exact sequence
$$0 \longrightarrow \Omega^{m-1}K \longrightarrow \Omega^{m-1}C\oplus \Omega^{m-1}E\oplus Q'_3 \longrightarrow \Omega^{m}X\oplus Q'_1 \longrightarrow 0,$$
\noindent where $Q'_i$ is projective, and $\Omega^{m-1}K$, $\Omega^{m-1}C$, $\Omega^{m-1}E$ are in finite representation type classes of modules. Then $\Lambda$ is $m$-Igusa-Todorov.
\end{dem}

\begin{remark} 	Note that \cite[Theorem 4.1 (3)]{W2} is a particular case of Theorem \ref{IT}. Indeed, let $\Lambda$ be an Artin algebra and $I,J$ be two ideals of $\Lambda$ such that $IJ\rad\,(\Lambda)=0.$ Assume that $\Lambda/I$ and $\Lambda/J$ are syzygy finite algebras with $I$ and 
$J$ of finite projective dimension as $\Lambda$-modules. Consider the following classes of 
$\Lambda$-modules $\mathcal{C}:=\{C\in \modu\,\Lambda\;|\; IC=0\}$, 
$\mathcal{D}:=\{D\in \modu\,\Lambda\;|\; JD=0\},$ $\mathcal{E}=\rm{add}\,\Lambda$ and $\mathcal{K}:=\Omega\mathcal{D}.$ It can be shown that the classes $\C,$ $\mathcal{D},$ $\mathcal{E}$ and $\mathcal{K}$ 
satisfy the needed hypothesis in Theorem \ref{IT}.
\end{remark}

\section{Lower triangular matrix Artin algebras}

In this section we study lower triangular matrix Artin algebras from the point of view of syzygy finite classes and the fact of being Igusa-Todorov. This kind of matrix algebras are of the form $\Lambda={\begin{pmatrix} T & 0 \\ {}_{U}M_T & U\end{pmatrix}},$ where $T$ and $U$ are Artin $R$-algebras,   $M$ is an $U$-$T$-bimodule, where $R$ acts centrally on $M,$ and furthermore  $M$ is 
finitely generated over $R.$ Note that, by \cite[Proposition III.2.1]{ARS},  it follows  that $\Lambda$ is also an Artin $R$-algebra. 
\

The description of the $\Lambda$-modules can be given, see the details in \cite[III.2]{ARS}, in terms of 
triples $(A,B,f),$ where $A$ is a $T$-module, $B$ is a $U$-module and $f:M\otimes A\to B$ is a morphism of $U$-modules. The pair $(\alpha,\beta):(A,B,f)\to(A',B',f')$ is a morphism of $\Lambda$-modules if 
$\alpha:A\to A'$ is a $T$-morphism, $\beta:B\to B'$ is a $U$-morphism and the following diagram commutes 
$$\xymatrix{ M\otimes A \ar[d]_f  \ar[r]^{M\otimes\alpha} & M\otimes A'\ar[d]^{f'}\\
B\ar[r]^\beta & B'.
}$$
A sequence $(A,B,f)\xrightarrow{(\alpha,\beta)}(A',B',f')\xrightarrow{(\alpha',\beta')}(A'',B'',f'')$ is exact if the sequences 
$A\xrightarrow{\alpha}A'\xrightarrow{\alpha'}A''$ and $B\xrightarrow{\beta}B'\xrightarrow{\beta'}B''$ are exact. Finally, the indecomposable 
projective $\Lambda$-modules are 
isomorphic to objects of the form $(P,M\otimes P,1)$ where $P$ is an indecomposable projective $T$-module, or of the form $(0,Q,0)$ where $Q$ is 
an indecomposable projective $U$-module \cite[Proposition III.2.5]{ARS}.
\

For a triangular matrix Artin algebra $\Lambda,$ as above, it is very useful to consider the following classes of $\Lambda$-modules 
$$\mathcal{T}:=\{\overline{A}:=(A,0,0)\in 
 \modu\,\Lambda\;:\; A\in \modu\,T\},$$ 
 $$\mathcal{U}:=\{\underline{B}:=(0,B,0)\in \modu\,\,\Lambda\;:\; B\in \modu\,U\}.$$ 
 As an applitacion of the Theorem \ref{IT}, we obtain the following result.
 
\begin{cor}\label{LT1}   Let $\Lambda={\begin{pmatrix} T & 0 \\ _{U}M_T & U\end{pmatrix}}$ be a triangular matrix Artin $R$-algebra such that $\T$ is $n$-syzygy finite and $\mathcal{U}$ is $m$-syzygy finite. Then, $\Lambda$ is a $k$-Igusa-Todorov algebra for $k:=\max\{m,n\}+1.$
\end{cor} 
\begin{dem} Consider the following classes of $\Lambda$-modules:  
$\mathcal{C}:=\mathcal{U}$, $\mathcal{D}:=\mathcal{T},$ $\mathcal{E}:=
\{(P, M\otimes P, 1)\;:\:  P\in \proj\,(T)\}\subseteq\proj\,(\Lambda)$ and $\mathcal{K}:=\Omega_\Lambda\mathcal{T}.$ Note that $\mathcal{U},$ $\mathcal{T}$ and $\mathcal{K}$ are $(k-1)$-syzygy finite 
classes in $\modu\,\Lambda.$ Moreover, the conditions (TS2) and (TS3) in Theorem \ref{IT} hold. 
\

For any $(A,B,f)\in\modu\,\Lambda,$ we have the exact sequence 
$$0\to \underline{B}\to (A,B,f)\to\overline{A}\to 0.$$
 Therefore $\modu\,\Lambda\subseteq \mathcal{U}*\mathcal{D}$ and thus the condition 
(TS1) in Theorem \ref{IT} holds. Hence $\Lambda$ is a $k$-Igusa-Todorov algebra.		
\end{dem}

In order to check the $n$-syzygy finiteness condition, for the classes $\mathcal{U}$ and $\mathcal{T},$ 
we need to have a description of $\Omega_\Lambda$ in terms of $\Omega_T$ and $\Omega_U.$ Under some 
natural hypothesis on the bimodule $M,$ we obtain the following useful description of $\Omega_\Lambda.$ To state this description, we 
consider the following notation. 
\

Let $\Gamma$ be an Artin algebra. For any $A\in\modu\,\Gamma,$ the exact sequence 
$$\cdots\to P^A_n\to P^A_{n-1}\to\cdots \to P^A_1\to P^A_0\to A\to 0$$ 
denotes the minimal projective resolution of $A.$ Thus, for any 
 $k\geq 1,$ we have an exact sequence $0\to\Omega^kA\xrightarrow{i^A_k} P^A_{k-1}\xrightarrow{\pi^A_{k-1}} \Omega^{k-1}A\to 0$ of $\Gamma$-modules.
 
\begin{lem}\label{triang}
Let $\Lambda={\begin{pmatrix} T & 0 \\ _{U}M_T & U\end{pmatrix}}$ be a triangular matrix Artin $R$-algebra such that $_UM \in\proj\,(U)$ and $M_T \in\proj\,(T).$ Then, the $n$-th syzygy of the $\Lambda$-module
$(A, B, f)$ is
{\small $$\Omega^n(A, B, f)=(\Omega^nA, M\otimes P^A_{n-1}, M\otimes i^A_{n})\oplus (0, \Omega^nB,0).$$}
\end{lem}
\begin{dem} Let $(A, B, f)\in\modu\,\Lambda.$ Then, the exact sequence 
$0\to\Omega A\xrightarrow{i^A_1} P^A_0\xrightarrow{\pi^A_0} A \to 0$ induces the following commutative diagram,
\[
\begin{tikzcd}[ampersand replacement=\&]
M\otimes \Omega A \arrow[dd] \arrow[rr,"M\otimes i^A_1", rightarrowtail] \&\& M\otimes P^A_0 \arrow[dd, "{\begin{pmatrix} 0\\ 1\end{pmatrix}}"] \arrow[rr, "{ M\otimes \pi^A_0 }", twoheadrightarrow] \&\& M\otimes A \arrow[dd, "f"]\\
\&\& \&\&\\
X \arrow[rr, "{\begin{pmatrix} \alpha\\ \beta \end{pmatrix}}"', rightarrowtail] \&\& P^B_0\oplus (M\otimes P^A_0) \arrow[rr,"{\begin{pmatrix} \pi^B_0 & f\circ (M\otimes \pi^A_0) \end{pmatrix}}"', twoheadrightarrow] \&\& B,\\
\end{tikzcd}
\]
\noindent where the first row is exact, since $M_T \in\proj\,(T),$ and $X$ is such that the second row is exact. Thus, from \cite[Corollary I.5.7]{ARS} it induces the pullback diagram

\[
\begin{tikzcd}[row sep=4ex, column sep=2ex]
0 \arrow[rr] && \Omega B \arrow[dd, equal] \arrow[rr] && X \arrow[dd, "-\alpha"] \arrow[rr, "\beta"] && M\otimes P^A_0 \arrow[dd, "f\circ (M\otimes \pi^A_0)"] \arrow[rr] && 0\\
&& && && &&\\
0 \arrow[rr] && \Omega B \arrow[rr] && P^B_0 \arrow[rr,"\pi^B_0"'] && B \arrow[rr] && 0 \\
\end{tikzcd}
\]
Since $_UM \in\proj\,(U)$ then $\Ext^1_U(M\otimes P^A_0,\Omega B)=0.$ Thus, the first row in the above diagram splits and then 
$X=\Omega B \oplus (M\otimes P^A_0).$ We get the following exact and commutative diagram
\[
\begin{tikzcd}[ampersand replacement=\&]
M\otimes \Omega A \arrow[dd,"{\begin{pmatrix} 0\\ M\otimes i^A_1\end{pmatrix}}"] \arrow[rr,"M\otimes i_A", rightarrowtail] \&\& M\otimes P^A_0 \arrow[dd, "{\begin{pmatrix} 0\\ 1\end{pmatrix}}"] \arrow[rr, "{ M\otimes\pi^A_0 }", twoheadrightarrow] \&\& M\otimes A \arrow[dd, "f"]\\
\&\& \&\& \\
\Omega B \oplus (M\otimes P^A_0)  \arrow[rr, "{\begin{pmatrix} \alpha_1 & \beta_1 \\ \alpha_2 &\beta_2 \end{pmatrix}}"', rightarrowtail] \&\& P_B\oplus (M\otimes P^A_0) \arrow[rr,"{\begin{pmatrix} 0 & f\circ (M\otimes\pi^A_0) \end{pmatrix}}"', twoheadrightarrow] \&\& B \\
\end{tikzcd}
\]
Thus, we get the equalities  
\begin{align*} \Omega(A, B, f) & =(\Omega A, \Omega B \oplus (M\otimes P^A_0), {\begin{pmatrix} 0\\ M\otimes i^A_1 \end{pmatrix}})\\
  & = (\Omega A, M\otimes P^A_0, M\otimes i^A_1)\oplus (0, \Omega B,0).
\end{align*}
Note that $_UM\otimes P_0^A \in\proj\,(U)$. In order to prove the result, we proceed by induction on $n.$ Indeed, suppose that 
{\small $$\Omega^n(A, B, f)=(\Omega^nA, M\otimes P^A_{n-1}, M\otimes i^A_{n}) \oplus (0, \Omega^nB,0).$$}
Then,  we have
{\small\begin{align*} \Omega^{n+1}(A, B, f) & =\Omega(\Omega^n A, M\otimes P^A_{n-1}, M\otimes i^A_n) \oplus(0,\Omega^{n+1}B,0).
\end{align*}}
On the other hand, by the equality $\Omega(X, Y, h)=(\Omega X, M\otimes P^X_0, M\otimes i^X_1)\oplus (0, \Omega Y,0),$ it follows
{\small\begin{align*}
\Omega(\Omega^n A, M\otimes P^A_{n-1}, M\otimes i^A_n) & =(\Omega^{n+1} A, M\otimes P^A_{n}, M\otimes i^A_{n+1}).
\end{align*}}
Finally, we conclude that 
{\small\begin{align*}
\Omega^{n+1}(A, B, f)&=(\Omega^{n+1}A, M\otimes P^A_{n}, M\otimes i^A_{n+1}) \oplus (0, \Omega^{n+1}B,0).
\end{align*} }
Thus, by induction we get the result.
\end{dem}

\begin{teo}\label{teotriang} Let $\Lambda={\begin{pmatrix} T & 0 \\ {}_{U}M_T & U\end{pmatrix}}$ be a triangular matrix Artin $R$-algebra such that $_UM \in\proj\,(U)$ and $M_T \in\proj\,(T).$ Then, 
for a positive integer $n,$ the following statements are 
equivalent.
\begin{itemize}
\item[(a)] The algebras $T$ and $U$ are $n$-syzygy finite.
\item[(b)] The classes $\mathcal{T}$ and $\mathcal{U}$ are $n$-syzygy finite.
\item[(c)] The algebra $\Lambda$ is $n$-syzygy finite.
\end{itemize}
Moreover, if one of the above equivalent conditions holds, then $\Lambda$ is a $(n-1)$-Igusa-Todorov algebra.
\end{teo}
\begin{dem} 
 Let $A\in\modu\,T$ and $B\in\modu\,U.$ Then, by Lemma \ref{triang}, we have
$\Omega^n_\Lambda\,\underline{B}=(0,\Omega^n_U B, 0)$ and 
$$\Omega^n_\Lambda\,\overline{A}=(\Omega^n_T A, M\otimes P^A_{n-1}, M\otimes i^A_n).$$
By the first equality, it is clear that $U$ is $n$-syzygy finite if, and only, if $\mathcal{U}$ is $n$-syzygy finite.
\

(a) $\Rightarrow$ (b) Let $C\in\modu\,T$ be such that $\Omega^n(\modu\,T)\subseteq\add\,(C).$  
Then, for any $A\in\modu\,T,$ we have that $M\otimes P^A_{k-1}\in
\add\,(M).$ Since $\Hom_U(M\otimes C,M)$ is an $R$-module of finite length, say $m,$ we fix a set 
$\{f_i\}_{i=1}^m$ of $R$-generators in $\Hom_U(M\otimes C,M).$ Then, for the family 
$\A:=\{(C,M,f_i)\}_{i=1}^m$ we get $\Omega^n_\Lambda(\mathcal{T})\subseteq \add\,(\A).$ Thus, 
the class $\mathcal{T}$ is $n$-syzygy finite.
\

(b) $\Rightarrow$ (a) Let $\C:=\{(C_i,D_i,f_i)\}_{i=1}^m$ be a family of indecomposable
$\Lambda$-modules such that $\Omega^n_\Lambda(\mathcal{T})\subseteq \add\,(\C).$ Let $A\in\modu\,T.$ 
Then $\Omega^n_\Lambda\,\overline{A}=\bigoplus_{i=1}^m\,(C_i,D_i,f_i)^{k_i}$ and thus
$\Omega^n_T\,A=\bigoplus_{i=1}^m\,C^{k_i}_i.$ Therefore $\Omega^n(\modu\,T)\subseteq \add\,(C),$ for 
$C:=\bigoplus_{i=1}^m\,C_i.$
\

(b) $\Leftrightarrow$ (c) It follows from the equality $\Omega^n(A,B,f)=\Omega^n\, \overline{A}\oplus \Omega^n\,\underline{B}.$\\ 

Finally, assume that $\Lambda$ is $n$-syzygy finite. Then, by \cite[Proposition 2.5]{W} we conclude that  $\Lambda$ is $(n-1)$-Igusa-Todorov.   
\end{dem}

Let $\Gamma$ be an Artin $R$-algebra. The $k$-triangular matrix $R$-algebra $T_k(\Gamma)$ is the matrix of 
size $k\times k$ with entries $[T_k(\Gamma)]_{i,j}:=\Gamma$ for $i\geq j,$ and 
$[T_k(\Gamma)]_{i,j}:=0$ otherwise.

\begin{cor}\label{Nsfin} For an Artin $R$-algebra $\Gamma,$ the following statements are equivalent.
\begin{itemize}
\item[(a)] $\Gamma$ is $n$-syzygy finite.
\item[(b)] $T_k(\Gamma)$ is $n$-syzygy finite, for any $k\geq 2.$
\item[(c)] $T_k(\Gamma)$ is $n$-syzygy finite, for some $k\geq 2.$
\end{itemize} 	
\end{cor}
\begin{dem}  (a) $\Rightarrow$ (b) Let $\Gamma$ be $n$-syzygy finite. We proceed by induction on $k.$ The case $k=2$ holds, by considering 
$M=\Gamma=T=U$ in Theorem \ref{teotriang}.
\

For the general case, we use that $T_{k+1}(\Gamma)={\begin{pmatrix} \Gamma & 0 \\ M & T_k(\Gamma)\end{pmatrix}}$, where $M$ is a matrix of size $k\times 1$ with entries 
$[M]_{i,1}:=\Gamma$ for any $i.$ Note that $M$ is a $T_k(\Gamma)-\Gamma$ bimodule, which is left and right projective. Assume by inductive hypothesis that $T_k(\Gamma)$ is $n$-syzygy finite. Then, by Theorem \ref{teotriang}, we get that $T_{k+1}(\Gamma)$ is n-syzygy finite. Thus, (b) is now true by induction.
\

(b) $\Rightarrow$ (c) It is trivial.
\

(c) $\Rightarrow$ (a) Apply Theorem \ref{teotriang} to the equality $T_{k}(\Gamma)={\begin{pmatrix} \Gamma & 0 \\ M & T_{k-1}(\Gamma)\end{pmatrix}}.$
\end{dem}

\begin{teo}\label{nITtr} Let $\Lambda={\begin{pmatrix} T & 0 \\ _{U}M_T & U\end{pmatrix}}$ be a triangular 
matrix Artin $R$-algebra such that $_UM \in\proj\,(U).$ and $M_T \in\proj\,(T).$ Then
$\Lambda$ is $n$-Igusa-Todorov if, and only, if $T$ and $U$ are $n$-Igusa-Todorov. 
\end{teo}
\begin{dem} $(\Rightarrow)$ Let $\Lambda$ be a $n$-Igusa-Todorov algebra, and let $(A,B,f)$ be a $n$-Igusa Todorov $\Lambda$-module. 
\

Let us prove that $U$ is a $n$-Igusa-Todorov algebra. Indeed, for $X \in\modu\,U,$ by Lemma \ref{ComIT} (b) there exist $\Lambda$-modules 
$(A_1,B_1,f_1)$, $(A_2,B_2,f_2)$ in $\add\,(A,B,f)$ and a short exact sequence 
$0 \to (A_2,B_2,f_2) \to (A_1,B_1,f_1) \to \Omega^n\,\underline{X} \to 0.$ Hence, we get the exact sequence $0 \to B_2 \to B_1 \to \Omega^n\,X \to 0$ with $B_1$, $B_2$ in $\add\,B.$ Then $U$ is a $n$-Igusa-Todorov algebra and $B$ is a $n$-Igusa-Todorov module.
\

We assert that $T$ is a $n$-Igusa-Todorov algebra. Indeed, let $Y \in\modu\,T.$ By Lemma \ref{triang} 
$\Omega^n_\Lambda\,\overline{Y}=(\Omega^n_T Y, M\otimes P^Y_{n-1}, M\otimes i^Y_{n}).$  
On the other hand, since $\Lambda$ is a $n$-Igusa-Todorov algebra, by Lemma \ref{ComIT} there is an exact sequence of $\Lambda$-modules
$0 \to (A'_2,B'_2,f'_2) \to (A'_1,B'_1,f'_1) \to \Omega^n_\Lambda\,\overline{Y}\to 0$
with $(A'_1,B'_1,f'_1)$ and $(A'_2,B'_2,f'_2)$ in $\add\,(A,B,f).$ Therefore, we get the exact sequence $0 \to A'_2 \to A'_1 \to \Omega^n_T Y \to 0$ of $T$-modules, with $A'_1$, $A'_2\in\add\,A;$ proving that $T$ is a $n$-Igusa-Todorov algebra and $A$ is a $n$-Igusa-Todorov module.
\

$(\Leftarrow)$ Assume that $T$ and $U$ are $n$-Igusa-Todorov algebras. Let $V\in\modu\,T$ and $V'\in\modu\,U$ be their fixed $n$-Igusa-Todorov modules. 
\
 
We prove that $\Lambda$ is a $n$-Igusa-Todorov algebra. Indeed, let $(A,B,f)\in\modu\,\Lambda.$ Hence we have an exact sequence $0 \to V_1 {\to} V_0 \xrightarrow{\rho} \Omega^n A \to 0$ of $T$-modules, and 
other exact sequence $0 \to V'_1 \to V'_0 \xrightarrow{\rho'} \Omega^n B \to 0$ of $U$-modules, with $V_i \in\add\,V$ and $V'_i \in \add\,V'.$ Then, we construct the following exact and commutative diagram of $U$-modules
 \[
 \begin{tikzcd}[ampersand replacement=\&]
 0 \arrow[dd]\& \& 0 \arrow[dd] \\
 \& \& \\
 M\otimes V_1 \arrow[rr, "\kappa"] \arrow[dd] \& \& V'_1 \arrow[dd]\\
 \& \& \\
 M\otimes V_0 \arrow[dd, "M\otimes \rho"] \arrow[rr,"{\begin{pmatrix} 0 \\ M\otimes i^A_n\rho  \end{pmatrix}}"]\& \& V'_0 \oplus (M\otimes P^A_{n-1}) \arrow[dd, "{\begin{pmatrix} \rho' & 0 \\ 0 & 1  \end{pmatrix}}"]\\
 \& \& \\
 M\otimes \Omega^n A \arrow[rr,"{\begin{pmatrix} 0 \\ M\otimes i^A_{n} \end{pmatrix}}"'] \arrow[dd]\& \&  \Omega^n B\oplus (M\otimes P^A_{n-1}) \arrow[dd]\\
 \& \& \\
 0\& \& 0.\\
 \end{tikzcd}
 \]
 Thus, we obtain the exact sequence of $\Lambda$-modules
 \
 
  $0 \to (V_1, V'_1, \kappa) \to (V_0, V'_0 \oplus (M\otimes P^A_{n-1}), {\begin{pmatrix} 0 \\ M\otimes i^A_{n}\rho \end{pmatrix}}) \to (\Omega^n A, \Omega^n B\oplus (M\otimes P^A_{n-1}), {\begin{pmatrix} 0 \\ M\otimes i^A_{n} \end{pmatrix}})\to 0$, which is the same as 
\ 

 $0 \to (V_1, V'_1, \kappa) \to (V_0, M\otimes P^A_{n-1}, M\otimes i^A_{n}\rho) \oplus (0, V'_0,0) \to (\Omega^n A, M\otimes P^A_{n-1}, M\otimes i^A_{n}) \oplus (0,\Omega^n B,0) \to 0.$ On the other hand, by 
 Lemma \ref{triang}
 {\small $$\Omega^n\,(A,B,f)=(\Omega^n A, M\otimes P^A_{n-1}, M\otimes i^A_{n}) \oplus (0,\Omega^n B,0).$$}
 Therefore, we obtain the exact sequence of $\Lambda$-modules
 $$0 \to (V_1, V'_1, \kappa) \to (V_0, M\otimes P^A_{n-1}, M\otimes i^A_{n}\rho) \oplus (0, V'_0,0)\to \Omega^n (A, B,f)\to 0.$$

\noindent Define $$C:=\displaystyle\bigoplus_{i=1}^r(V,V'\oplus M,h_i),$$
 where $h_1,h_2,\cdots,h_r$ are $R$-generators of $\Hom_U(M\otimes V,V'\oplus M).$ 
 It follows that $C$ is a $n$-Igusa-Todorov $\Lambda-$module.
\end{dem}

\begin{cor}\label{LT2} For an Artin $R$-algebra $\Gamma,$ the following statements are equivalent. 
\begin{itemize}
\item[(a)] $\Gamma$ is $n$-Igusa-Todorov.
\item[(b)] $T_k(\Gamma)$ is $n$-Igusa-Todorov, for any $k\geq 2.$
\item[(c)] $T_k(\Gamma)$ is $n$-Igusa-Todorov, for some $k\geq 2.$
\end{itemize}
\end{cor}
\begin{dem} By using Theorem \ref{nITtr}, we can proceed in a very similar way as we did in the proof of 
Corollary \ref{Nsfin}.
\end{dem}

\bibliographystyle{unsrt}

\end{document}